\documentclass[runningheads,a4paper]{llncs}

\usepackage{amssymb}
\usepackage{graphicx}
\usepackage[applemac]{inputenc}
\usepackage{amsmath}
\usepackage{color}
\usepackage{subfigure}  

\usepackage{url}
\newcommand{\keywords}[1]{\par\addvspace\baselineskip
\noindent\keywordname\enspace\ignorespaces#1}

\def\be{\begin{equation}}
\def\ee{\end{equation}}
\def\bea{\begin{eqnarray}}
\def\eea{\end{eqnarray}}

\begin{document}

\mainmatter  

\title{On the optimal control of opinion dynamics on evolving networks}
\titlerunning{On the optimal control of opinion dynamics on evolving networks}

\author{Giacomo Albi\inst{1} \and Lorenzo Pareschi\inst{2} \and Mattia Zanella\inst{2}}
\authorrunning{G. Albi, L. Pareschi, M. Zanella}
\institute{
TU M\"unchen, Faculty of Mathematics,\\
Boltzmanstra\ss e 3, D-85748, Garching (M\"unchen), Germany.
\and
Department of Mathematics and Computer Science, \\
Via N. Machiavelli 35, 44121, University of Ferrara, Italy.
}
\maketitle

\begin{abstract}
In this work we are interested in the modelling and control of opinion dynamics spreading on a time evolving network with scale-free asymptotic degree distribution. The mathematical model is formulated as a coupling of an opinion alignment system with a probabilistic description of the network. The optimal control problem aims at forcing consensus over the network, to this goal a control strategy based on the degree of connection of each agent has been designed. A numerical method based on a model predictive strategy is then developed and different numerical tests are reported. The results show that in this way it is possible to drive the overall opinion toward a desired state even if we control only a suitable fraction of the nodes.
\end{abstract}

\section{Introduction}
Graph theory has emerged in recent years as one of the most active fields of research \cite{AB,BAJ,BBV,XZW}. In fact, the study of technological and communication networks earned a special attention thanks to a huge amount of data coming from empirical observations and more recently from online platforms like Facebook, Twitter, Instagram and many others. This fact offered a real laboratory for testing on a large-scale the collective behavior of large populations of agents \cite{KGH,S} and new challenges for the scientific research has emerged. In particular, the necessity to handle millions, and often billions, of vertices implied a substantial shift to large-scale statistical properties of graphs giving rise to the study of the so-called scale-free networks \cite{BAJ,N,XZW}. 

In this work, we will focus our attention on the modelling and control of opinion dynamics on a time evolving network.
We consider a system of agents, each one belonging to a node of the network, interacting only if they are connected through the network. Each agent modifies his/her opinion through a compromise function which depends both on opinions and the network \cite{AHP,APZa,APZb,C,DGM,W}. At the same time new connections are created and removed from the network following a preferential attachment process. For simplicity here we restrict to non-growing network, that is a graph where the total number of nodes and the total number of edges are conserved in time. An optimal control problem is then introduced in order to drive the agents toward a desired opinion. 
The rest of the paper is organized as follows. In Section \ref{sec:micro} we describe the alignment model for opinions spreading on a non-growing network. In order to control the trajectories of the model we introduce in Section \ref{sec:control} a general setting for a control technique weighted by a function on the number of connections. A numerical method based on model predictive control is then developed. Finally in Section \ref{sec:numerics} we perform numerical experiments showing the effectiveness of the present approach. Some conclusion are then reported in the last Section. 
\section{Modelling opinion dynamics on networks}\label{sec:micro}
In the model each agent $i=1,\dots,N$ is characterized by two quantities $(w_i,c_i), i=1,\dots,N$, representing the opinion and the number of connections of the agent $i$th respectively. This latter term is strictly related to the architecture of the social graph where each agent shares its opinion and influences the interaction between individuals. Each agent is seen here as a node of a time evolving graph $\mathcal{G}^N=\mathcal{G}^N(t), t\in[t_0,t_f]$ whose nodes are connected through a given set of edges. In the following we will indicates the density of connectivity the constant $\gamma\ge 0$.
\subsection{Network evolution without nodes' growth}\label{sec:net_growth}
In the sequel we will consider a graph with both a fixed number of nodes $N$ and a fixed number of edges $E$. In order to describe the network's evolution we take into account a preferential attachment probabilistic process. This mechanism, known also as  Yule process or Matthew effect, has been used in the modeling of several phenomena in biology, economics and sociology, and it is strictly connected to the generation of power law distributions \cite{BAJ,XZW}. The initial state of the network, $\mathcal{G}^N(0)$, is chosen randomly and, at each time step an edge is randomly selected and removed from the network. At the same time, a node is selected with probability
\begin{equation}\label{eq:probability_pi}
\Pi_{\alpha}(c_i) = \dfrac{c_i+\alpha}{\sum_{j=0}^N (c_j+\alpha)}=\dfrac{c_i+\alpha}{2E+N\alpha},\qquad i=1,\dots,N,
\end{equation}
among all possible nodes of $\mathcal{G}^N$, with $\alpha>0$ an attraction coefficient. Based on the  probability \eqref{eq:probability_pi} another node is chosen at time $t$ and connected with the formerly selected one. The described process is repeated at each time step. In this way both the number of nodes and the total number of edges remains constant in the reference time interval. Let $p(c,t)$ indicates the probability that a node is endowed of degree $c$ at time $t$, we have
\begin{equation}\begin{split}
\sum_c p(c,t)=1, \qquad
\sum_c c~p(c,t)=\gamma.
\end{split}\end{equation}
Then we have that the described process is in agreement with the following master equation 
\begin{equation}\begin{split}\label{eq:master}
\dfrac{d}{dt} p(c,t)=&\dfrac{D}{E}\left[(c+1)p(c+1,t)-cp(c,t)\right]\\
&+\dfrac{2D}{2E+N\alpha}\left[(c-1+\alpha)p(c-1,t)-(c+\alpha)p(c,t)\right],
\end{split}\end{equation}
where $D>0$ characterizes the relaxation velocity of the network toward an asymptotic degree distribution $p_{\infty}(c)$, the righthand side consists of four terms, the first and the third terms account the rate of gaining a node of degree $c$ and respectively the second and fourth terms the rate of losing a node of degree $c$. The equation \eqref{eq:master} holds in the interval $c\le E$, whereas for each $c>E$ we set $p(c,t)=0$. While most the random graphs models with fixed number of nodes and vertices produces unrealistic degree distributions as the Watts and Strogatz generation model, called small-world model \cite{WS}, the main advantage of the graph generated through the described process is the possibility to recover the scale-free properties. Indeed we can easily show that if $\gamma=2E/N \ge 1$ with attraction coefficient $\alpha\ll 1$ then the stationary degree distribution $p_{\infty}(c)$ obeys a power-law of the following form
\begin{equation}
p_{\infty}(c)=\left(\dfrac{\alpha}{\gamma}\right)^{\alpha}\dfrac{\alpha}{c}.
\end{equation}
When $\alpha\gg 1$ we loose the features of the preferential attachment mechanism, in fact high degree nodes are selected approximately with the same probability of the nodes with low degree of connection. Then the selection occurs in a non preferential way and the asymptotic degree distribution obeys the Poisson distribution
\begin{equation}
p_{\infty}(c)=\dfrac{e^{-\gamma}}{c!}\gamma^c.
\end{equation}

\begin{figure}
\centering
\includegraphics[scale=0.2]{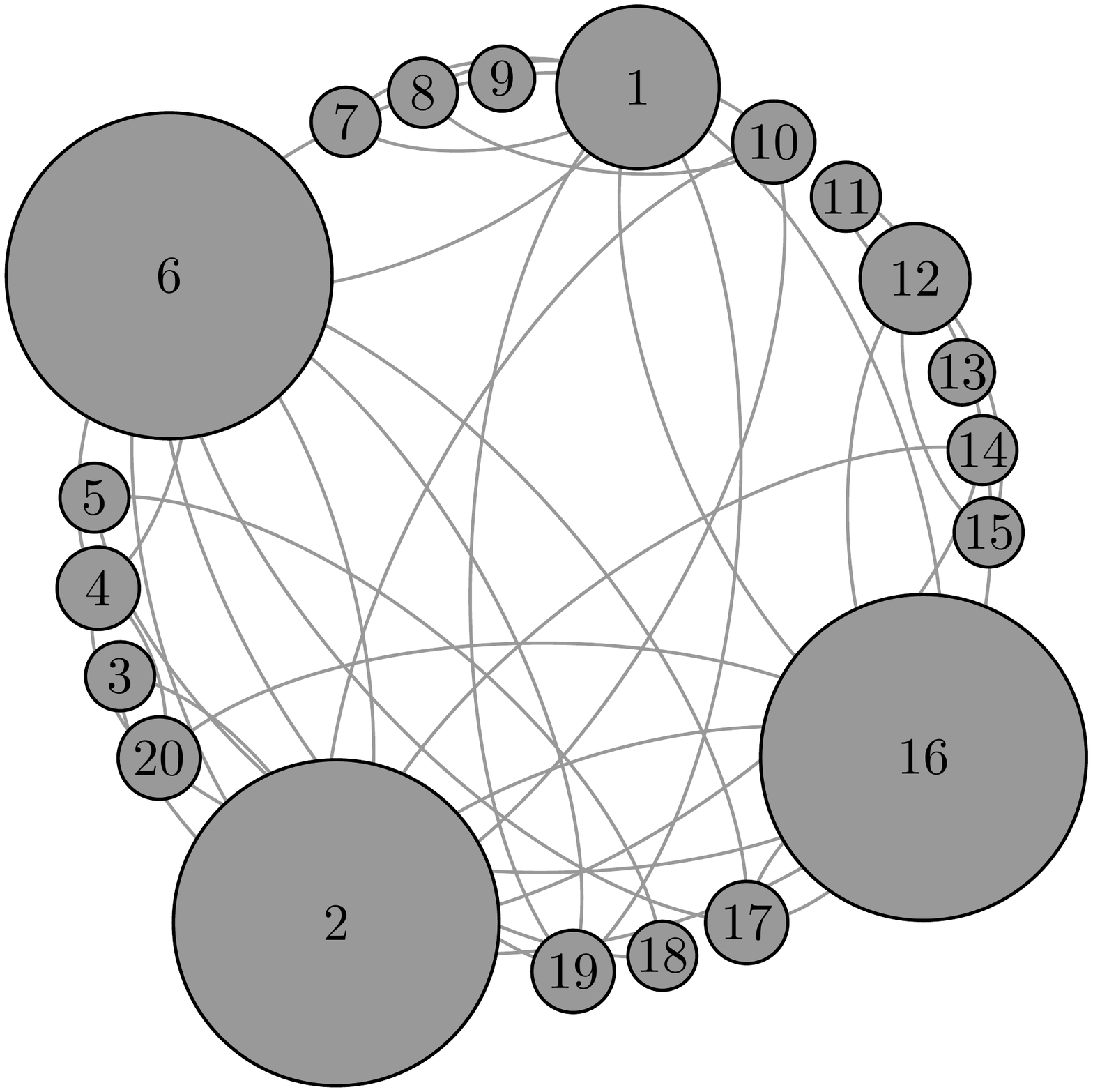}
\hspace{+1.5cm}
\includegraphics[scale=0.2]{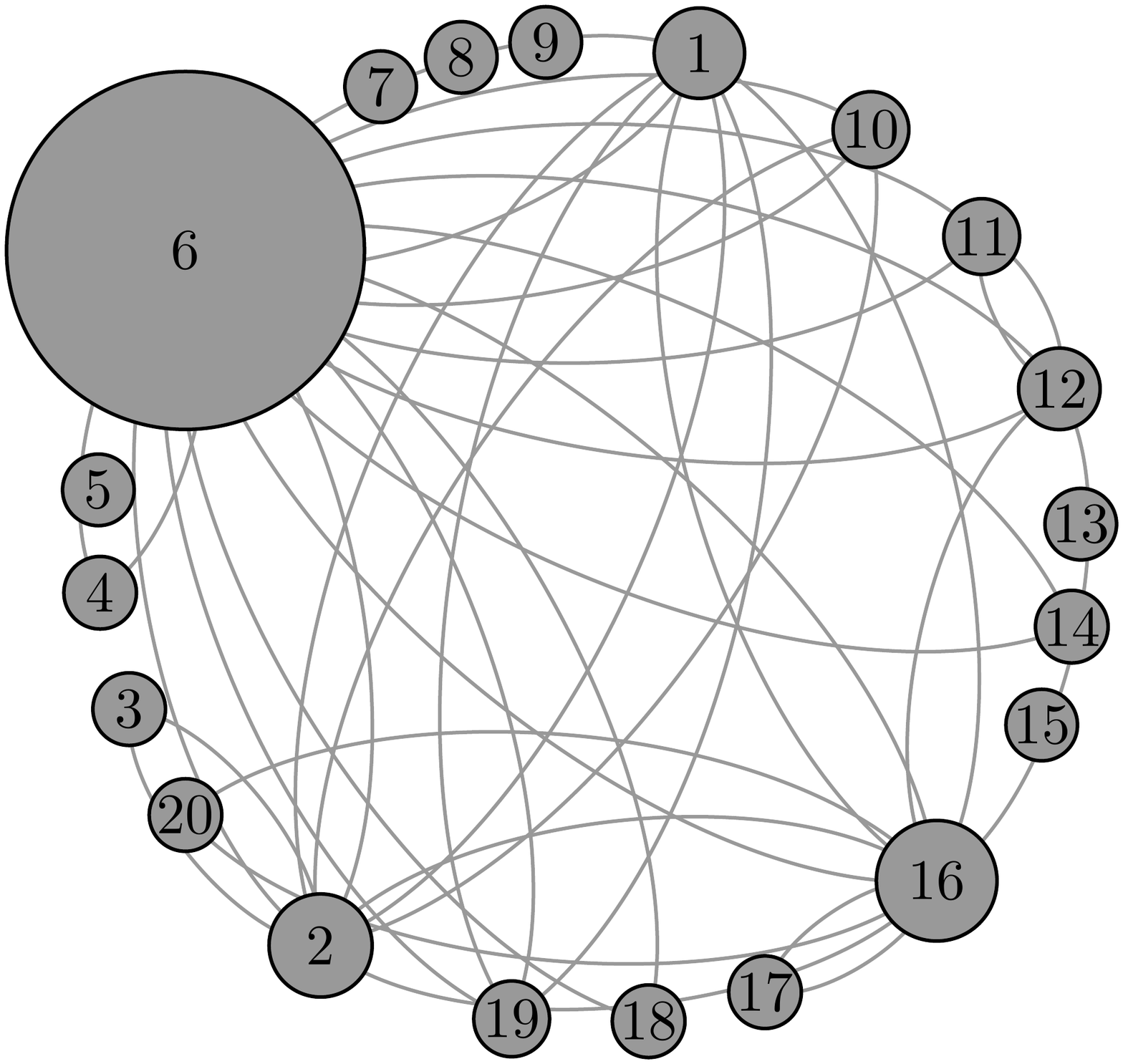}
\caption{Left: initial configuration of the sample network $\mathcal{G}^{20}$ with density of connectivity  $\gamma=5$. Right: a simulation of the network $\mathcal{G}^{20}$ after $10$ time steps of the preferential attachment process. The diameter of each node is proportional to its degree of connection.}
\label{fig:evo_graph_1}
\end{figure}

A simple graph is sketched in Figure \ref{fig:evo_graph_1} where we can observe how the initial degree of the nodes influences the evolution of the connections. In order to correctly observe the creation of the new links, that preferentially connect nodes with the highest connection degree, we marked each node with a number $i=1,\dots,20$ and the nodes' diameters are proportional with their number of connections. 

\subsection{The opinion alignment dynamics}
The opinion of the $i$th agent ranges in the closed set $I=[-1,1]$, that is $w_i=w_i(t)\in I$ for each $t\in [t_0,t_f]$, and its opinion changes over time according to the following differential system
\begin{equation}\label{eq:dynamics}
\dot{w}_i = \frac{1}{|S_i|} \sum_{j\in S_i} P_{ij} (w_j-w_i), \qquad i=1,\dots,N
\end{equation}
where $S_i$ indicates the set of vertex connected with the $i$th agent and reflects the architecture of the chosen network, whereas $c_i=|S_i| < N$ stands for the cardinality of the set $S_i$, also known as degree of vertex $i$. Note that the number of connections $c_i$ evolves in time accordingly to the process described in Section \ref{sec:net_growth}. Furthermore we introduced the interaction function $P_{ij}\in[0,1]$, depending on the opinions of the agents and the graph $\mathcal{G}^N$ which can be written as follows
\begin{equation}
P_{ij}=P(w_i,w_j; \mathcal{G}^N).
\end{equation}
A possible choice for the interaction function is the following
\begin{equation}\label{eq:interaction_fun}
P(w_i,w_j; \mathcal{G}^N)=H(w_i,w_j)K(\mathcal{G}^N),
\end{equation}
where $H(\cdot,\cdot)$ represents the  positive compromise propensity, 
and $K$ a  general function taking into account statistical properties of the graph $\mathcal{G}$. In what follows we will consider $K=K(c_i,c_j)$, a function depending on the vertices’ connections.

\section{Optimal control problem of the alignment model}\label{sec:control}
In this section we introduce a control strategy which characterizes the action of an external agent with the aim of driving opinions toward a given target $w_d$. To this goal, we consider the evolution of the network $\mathcal{G}^N(t)$ and the opinion dynamics in the interval $[t_0,t_f]$. Therefore we introduce the following optimal control problem
\begin{equation}\label{eq:functional}
\min_{u\in\mathcal{U}}J(\textbf{w},u):= \frac{1}{2}\int_{t_0}^{t_f} \Big\{ \frac{1}{N}\sum_{j=1}^N (w_j(s)-w_d)^2 +\nu u(s)^2  \Big\} ds,
\end{equation}
subject to 
\begin{equation}\begin{split}\label{eq:constrain}
\dot w_i &= \dfrac{1}{|S_i|}\sum_{j\in S_i}P_{ij}(w_j-w_i)+u\chi(c_i\ge c^*),\quad
w_i(0)=w_i^0,
\end{split}\end{equation}
where we indicated with $\mathcal{U}$ the set of admissible controls, with $\nu>0$ a regularization parameter which expresses the strength of the control in the overall dynamics and $w_d\in[-1,1]$ the target opinion. 
Note that the action of the control $u$ is weighted by an indicator function $\chi(\cdot)$, which is active only for  the nodes with degree $c_i\ge c^*$.
In general this selective control approach models an a-priori strategy of a policy maker, possibly acting under limited resources or unable to influence the whole ensemble of agents.

The solution of this kind of control problems is in general a difficult task, given that their direct solution is prohibitively expensive for a large number of agents. Different strategies have been developed for alignment modeling in order to obtain feedback controls or more general numerical control techniques \cite{ABCK,AHP,APZa,APZb,BF,WCB}. 
To tackle numerically the described problem a standard strategy makes use of a model predictive control (MPC) approach, also referred as receding horizon strategy.

In general MPC strategies solves a finite horizon open-loop optimal control problem predicting the dynamic behavior over a predict horizon $t_p\leq t_f$, with initial state sampled at time $t$ (initially $t = t_0$), and computing the control on a control horizon $t_c\leq t_p$. The optimization is computed introducing
a new integral functional $J_p(\cdot,\cdot)$, which is an approximation of \eqref{eq:functional} on the time interval $[t,t+t_p]$, namely
\begin{equation}\label{eq:sub_functional}
J_p(\textbf{w},\bar{u}):= \frac{1}{2}\int_t^{t+t_p} \Big\{ \frac{1}{N}\sum_{j=1}^N (w_j(s)-w_d)^2 +\nu_p \bar{u}(s)^2  \Big\} ds
\end{equation}
where the control, $\bar{u}: [t,t+t_p]\to \mathcal{U}$, is supposed to be an admissible control in the set of admissible control $\mathcal{U}$, subset of $\mathbb{R}$, and $\nu_p$ a possibly different penalization parameter with respect to the full optimal control problem.
Thus the computed optimal open-loop control $\bar{u}(\cdot)$ is applied feedback to the system dynamic until the next sampling time $t+t_s$ is evaluated, with $t_s\leq t_c$, thereafter the procedure is repeated taking as initial state of the dynamic at time $t+t_s$ and shifting forward the prediction and control horizons, until the final time $t_f$ is reached. This process generates a sub-optimal solution with respect to the solution of the full optimal control problem \eqref{eq:functional}-\eqref{eq:constrain}.

Let us consider now the full discretize problem,  defining the time sequence $[t_0,t_1,\ldots,t_M]$, where  $t_{n}-t_{n-1}=t_s=\Delta t>0$ and $t_M:=M\Delta t = t_f$, for all $n=1,\ldots, M $, assuming furthermore that $t_c = t_p = p\Delta t$, with $p>0$. Hence the linear MPC method look for a piecewise control on the  time frame $[t_0,t_M]$, defined as follows
\begin{align}
\bar{u}(t) = \sum_{n=0}^{M-1} \bar{u}^n\chi_{[t_n,t_{n+1}]}(t).
\end{align} 
In order to discretize the evolution dynamic we consider a Runge-Kutta  scheme, the full discretized optimal control problem on the time frame $[t_n,t_n+p\Delta t]$ reads 
\begin{equation}\label{eq:disc_functional}
\min_{\bar{u}\in\mathcal{U}}J_p(\textbf{w},\bar{u}):= \frac{1}{2}\int_{t_n}^{t_n+p\Delta t} \Big\{ \frac{1}{N}\sum_{j=1}^N (w_j(s)-w_d)^2 +\nu_p \bar{u}^2  \Big\} ds
\end{equation}
subject to 
\begin{equation}
\begin{subequations}
\begin{aligned}
W^{(n)}_{i,l} &= w_{i}^n   +  \Delta t\sum_{k=1}^s a_{l,k} \left(F(t+\theta_k\Delta t,W^{(n)}_{i,k})+\bar{U}^{(n)}_kQ_i(t+\theta_k\Delta t)\right),\\
w_i^{n+1} &= w_{i}^n +  \Delta t\sum_{l=1}^s b_{l} \left(F(t+\theta_l\Delta t,W^{(n)}_{i,l})+\bar{U}^{(n)}_lQ_i(t+\theta_l\Delta t)\right),\\
w_i^n &= w_i(t_n),
\end{aligned}
\end{subequations}
\end{equation}
for all $n=1,\ldots, p-1$; $l=1,\ldots,s$; $i,\ldots, N$ and having defined the following functions
\begin{align*}
F(t,w_i)=\dfrac{1}{|S_i(t)|}\sum_{j\in S_i(t)}P_{ij}(w_j-w_i),\quad
Q_i(t) = \chi(c_i(t)\ge c^*).
\end{align*}
The coefficients  $(a_{l,k})_{l,k}$, $(b_l)_l$ and $(\theta_{l})_l$, with $l,k=1,\ldots,s$, define the Runge-Kutta method and $(\bar{U}^{(n)})_l, (W^{(n)}_{i,l})_l$ are the internal stages associated to $\bar{u}(t), w_i(t)$ on time frame $[t_n,t_{n+1}]$.

\begin{remark}{(Instantaneous control).}\label{rem:inst}
Let us restrict to the case of a single prediction horizon, $p = 1$, where we discretize the dynamic with an explicit Euler scheme ( $a_{1,1}=\theta_1=0$ and $ b_1 = 1$). Notice that since the control $\bar{u}$ is a constant value and assuming that the network, $\mathcal{G}^N$ remains fixed over the time interval $[t_n,t_n+\Delta t]$ the discrete optimal control problem \eqref{eq:disc_functional}  reduces to
\begin{equation}
\min_{\bar{u}\in\mathcal{U}}J_p(\textbf{w},\bar{u}^n):= \Delta t \Big\{ \frac{1}{N}\sum_{j=1}^N (w_j^{n+1}(\bar{u}^n)-w_d)^2 +\nu_p (\bar{u}^n)^2  \Big\}
\end{equation}
with 
\begin{align}\label{eq:EE}
w_i^{n+1} &= w_{i}^n +  \Delta t\left(F(t_n,w^n_i)+\bar{u}^nQ_i^n\right),\quad w_i^n= w_i(t_n).
\end{align}
In order to find the minima of \eqref{eq:disc_functional} is sufficient to find the value $\bar{u}$ satisfying $ \partial_{\bar{u}} J_p(\textbf{w},\bar{u})=0$, which can be computed by  a straightforward calculation
\begin{align}\label{eq:IC}
\bar{u}^n = -\frac{1}{N\nu+ \Delta t\sum_{j=1}^N(Q^n_j)^2}\left( \sum_{j=1}^NQ^n_j\left(w^n_j-w_d\right)+\Delta t\sum_{j=1}^NQ^n_jF(t_n,w^n_i)\right).
\end{align}
where we scaled the penalization parameter with $\nu_p= \Delta t \nu$. 
\end{remark}

\section{Numerical results}\label{sec:numerics}
\begin{figure}
\centering
\subfigure[$c^*=10$]
{\includegraphics[scale=0.22]{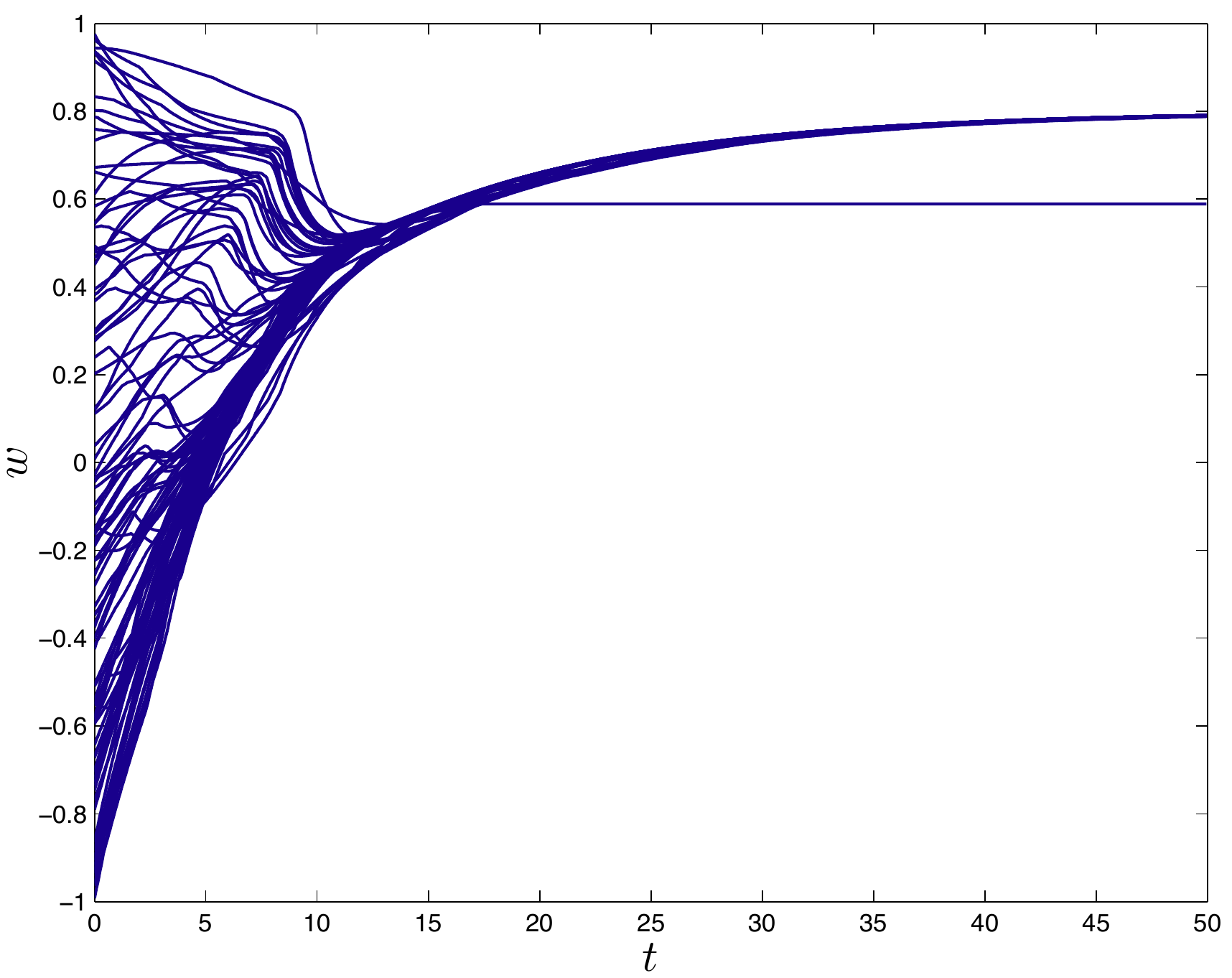}}
\subfigure[$c^*=20$]
{\includegraphics[scale=0.22]{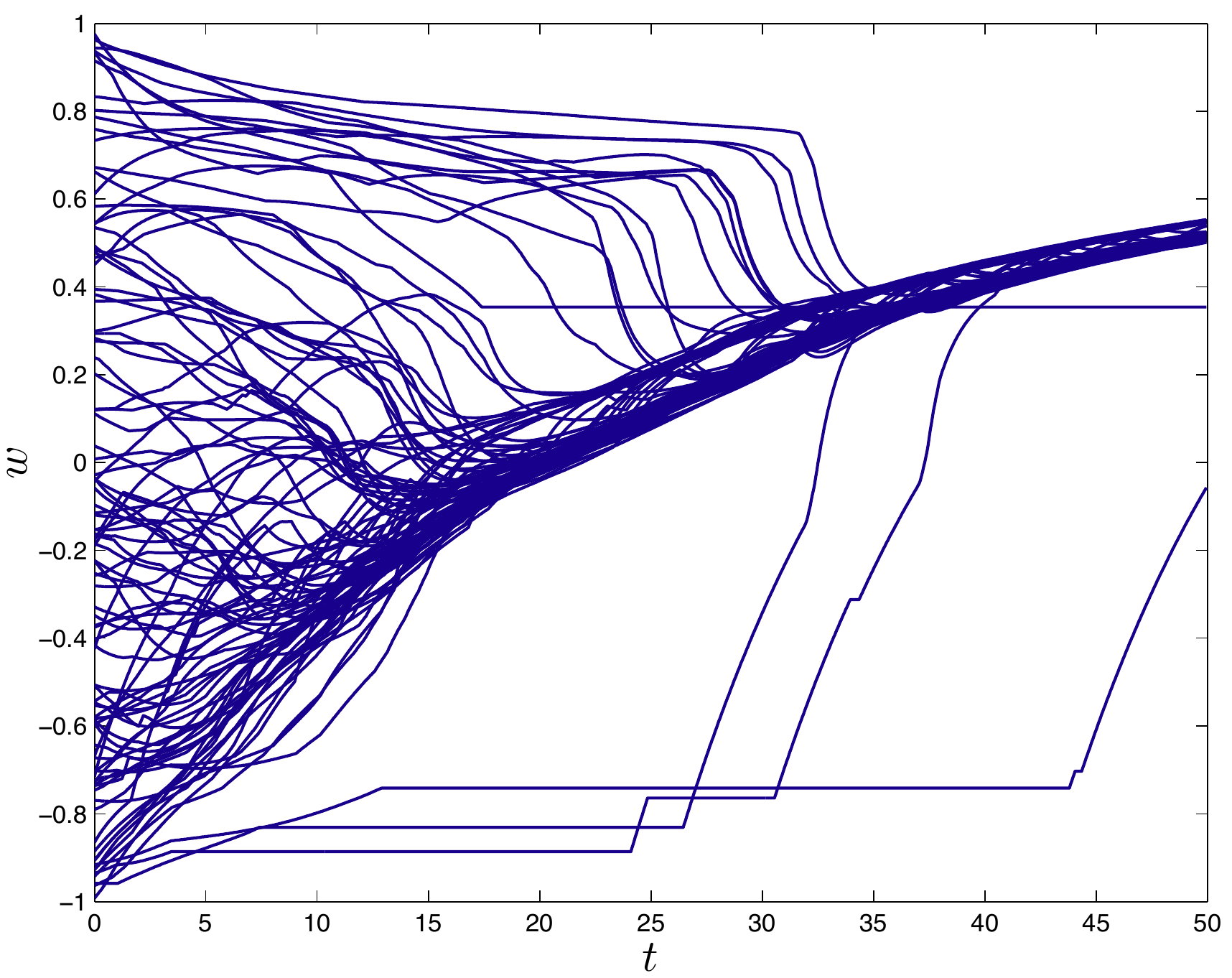}}
\subfigure[$c^*=30$]
{\includegraphics[scale=0.22]{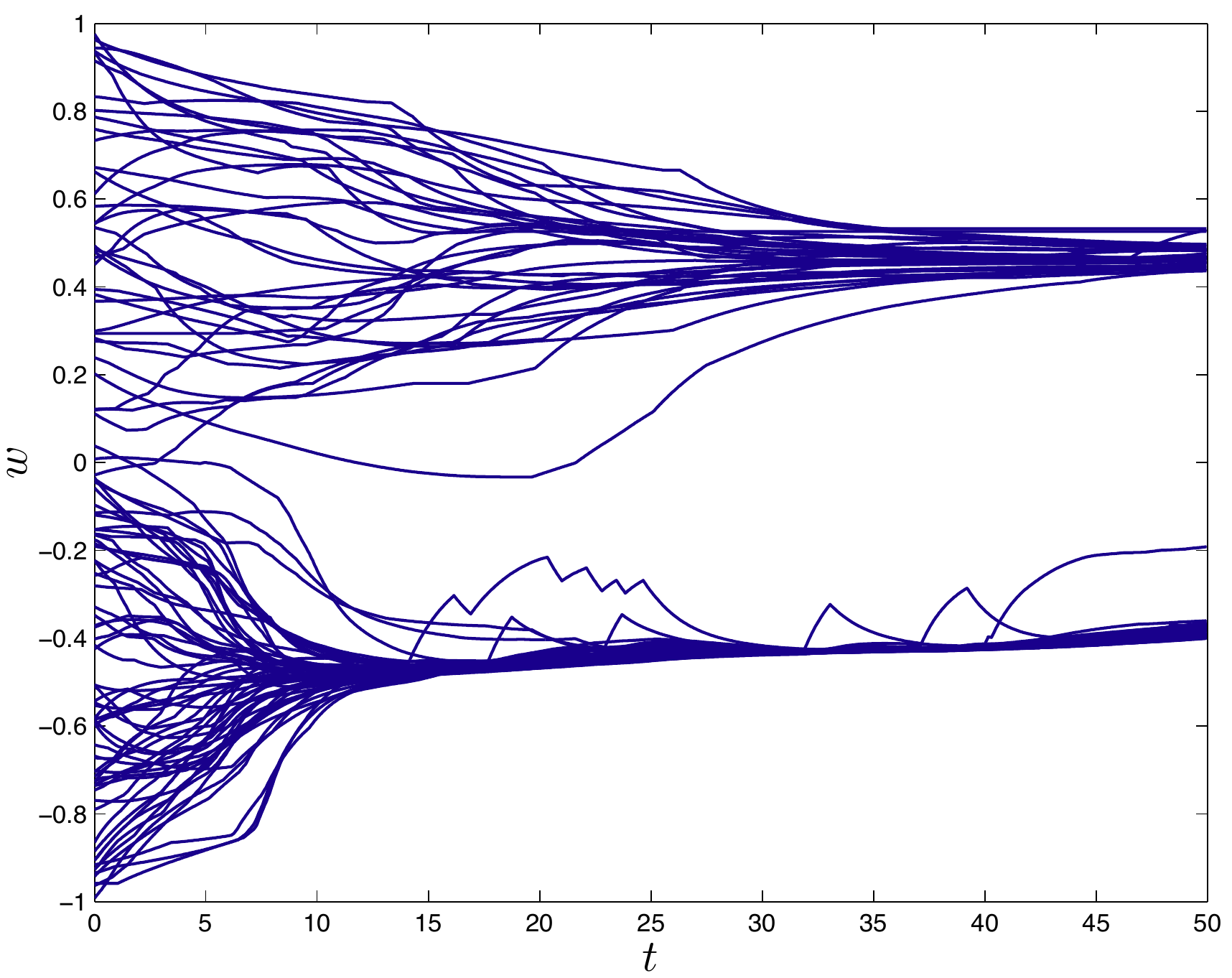}}
\caption{Evolution of the constrained opinion dynamics with uniform initial distribution of opinions over the time interval $[0, 50]$ for different values of $c^*=10,15,30$ with target opinion $w_d=0.8$, control parameter $\kappa=0.1$, $\Delta t=10^{-3}$ and confidence bound $\Delta=0.4$.} \label{fig:line_nonlin2}
\end{figure}

\begin{figure}
\centering
\includegraphics[scale=0.25]{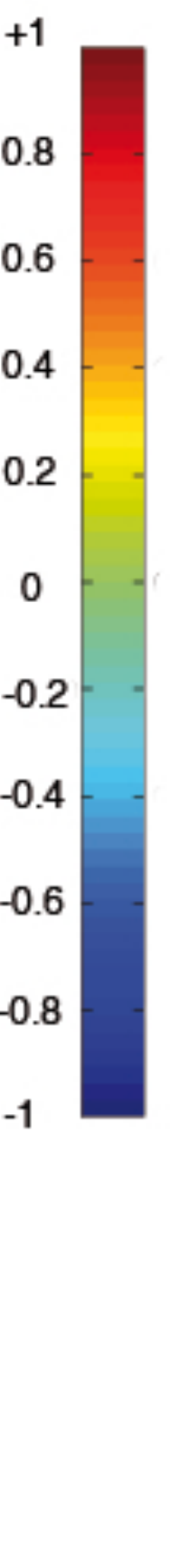}
\includegraphics[scale=0.18]{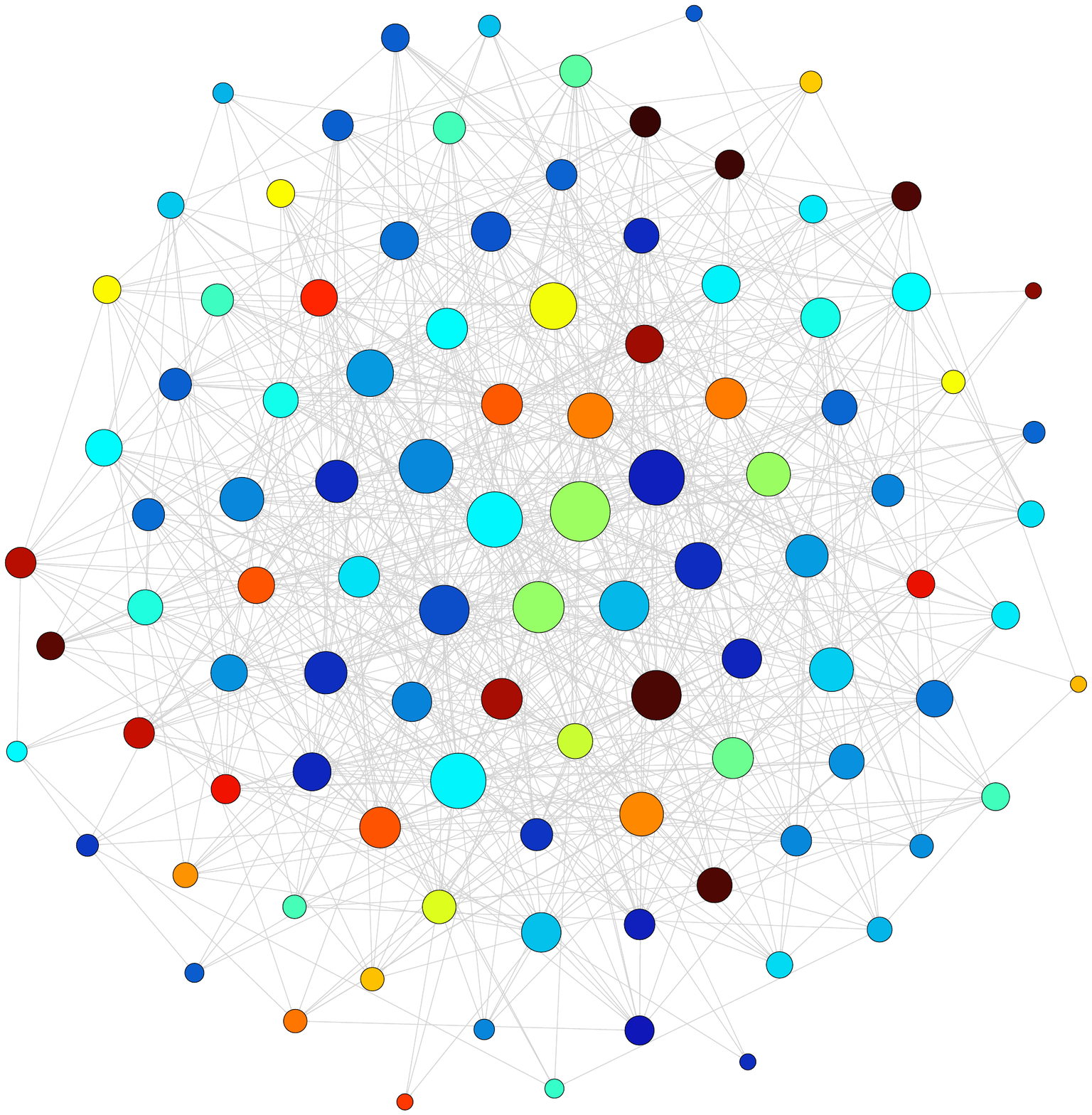}
\includegraphics[scale=0.18]{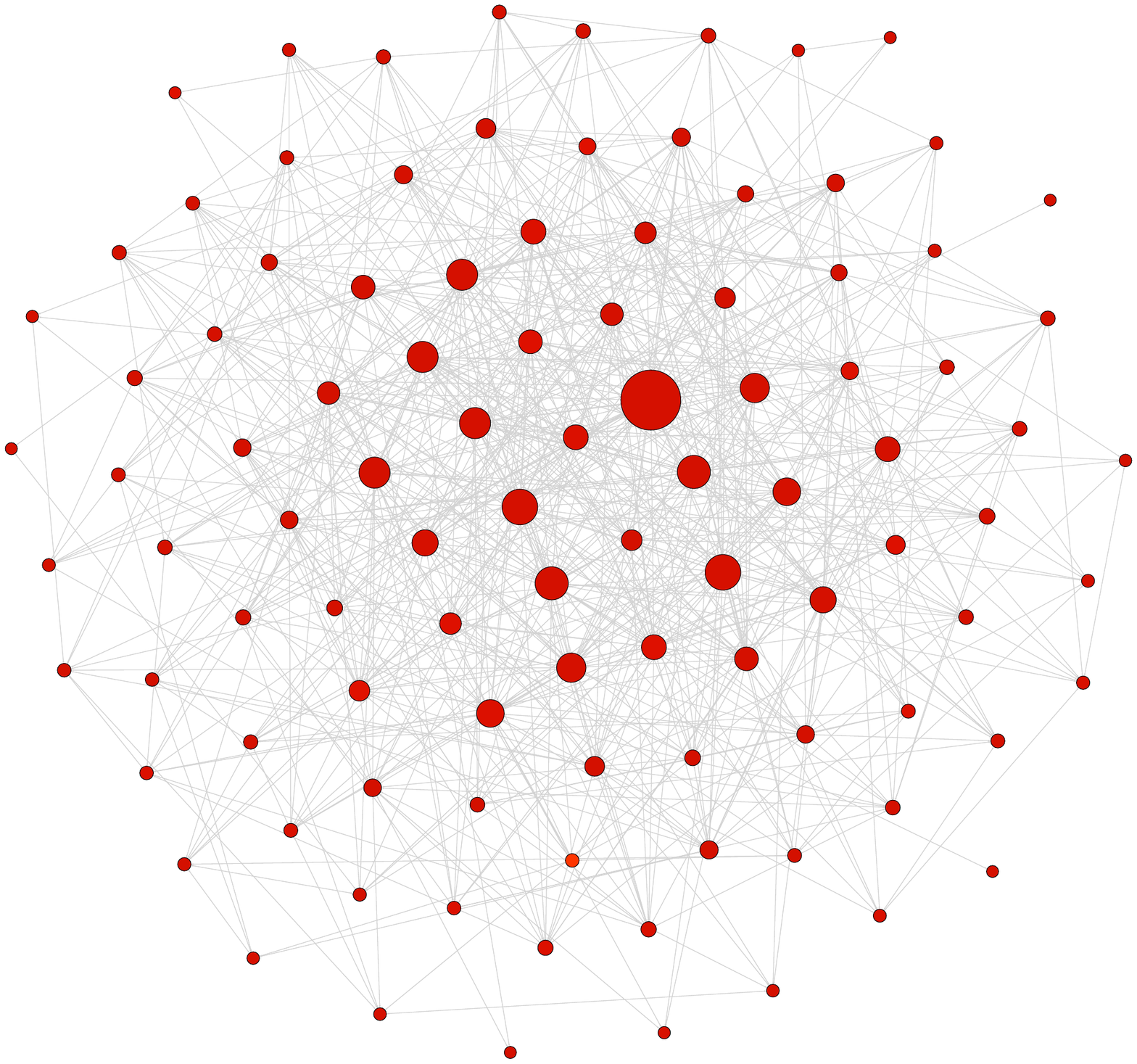}
\includegraphics[scale=0.18]{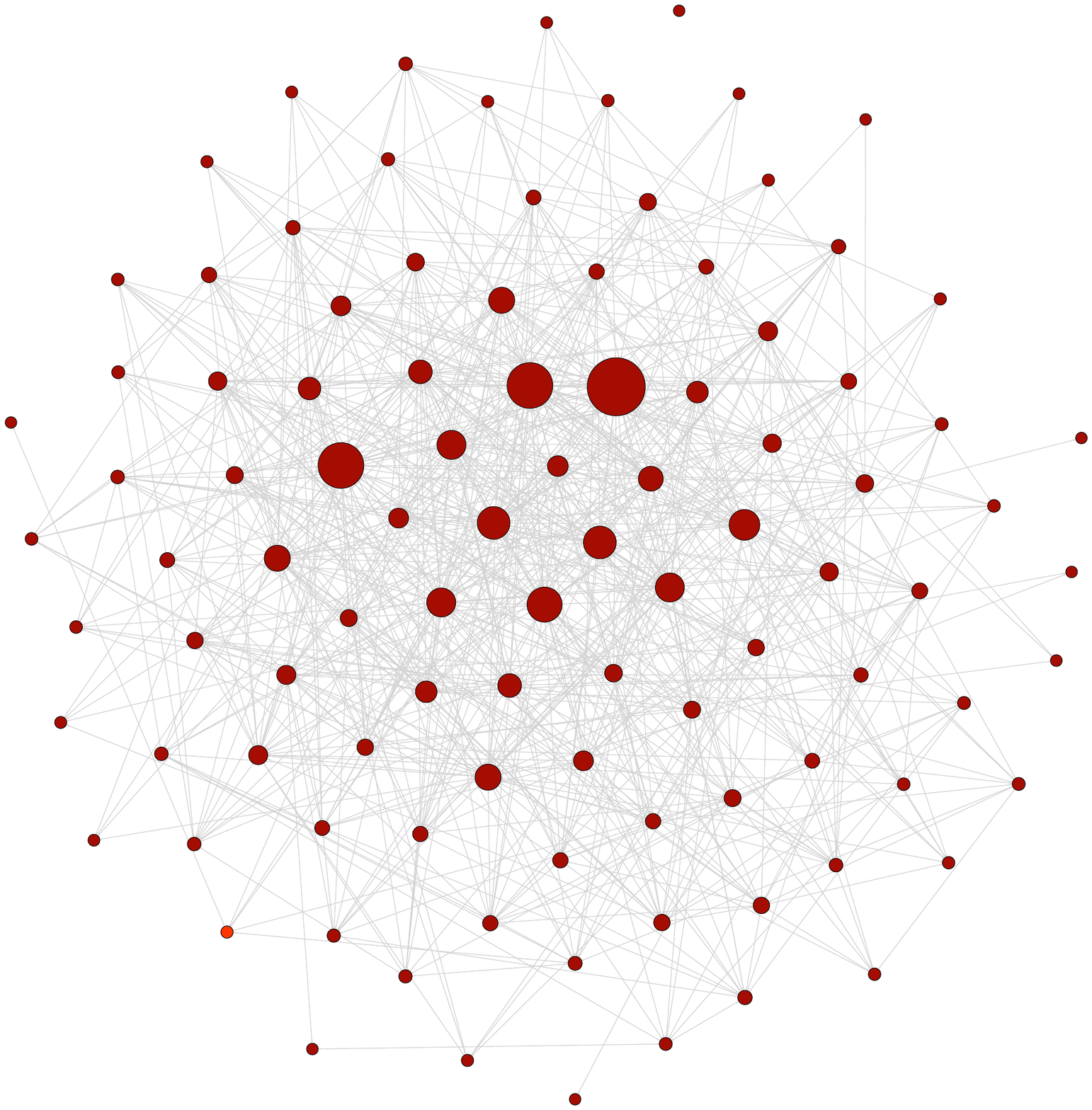} \\
\includegraphics[scale=0.25]{colorbar_modified}
\includegraphics[scale=0.18]{t0_Knonlin}
\includegraphics[scale=0.18]{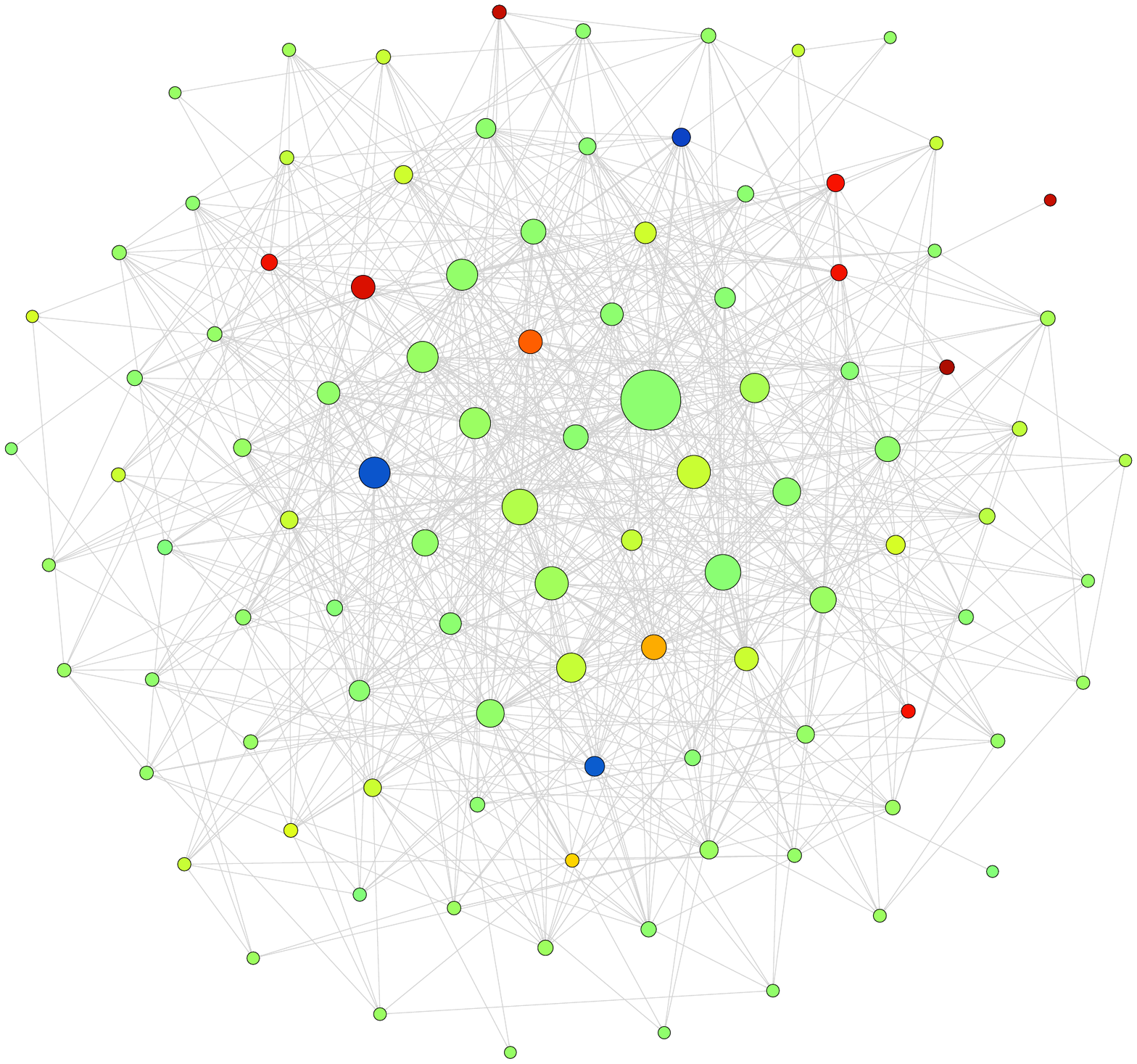}
\includegraphics[scale=0.18]{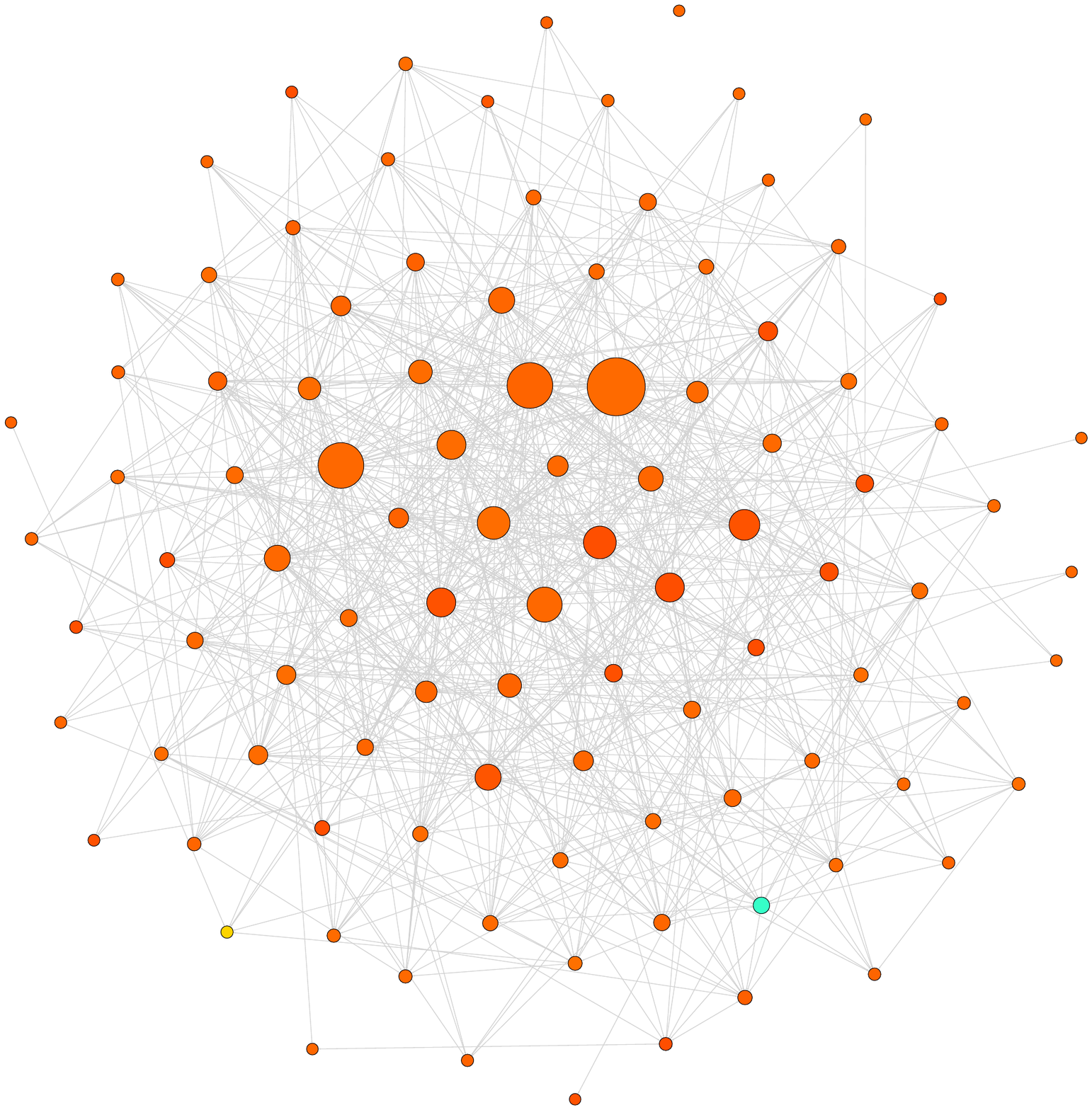} \\
\includegraphics[scale=0.25]{colorbar_modified}
\includegraphics[scale=0.18]{t0_Knonlin}
\includegraphics[scale=0.18]{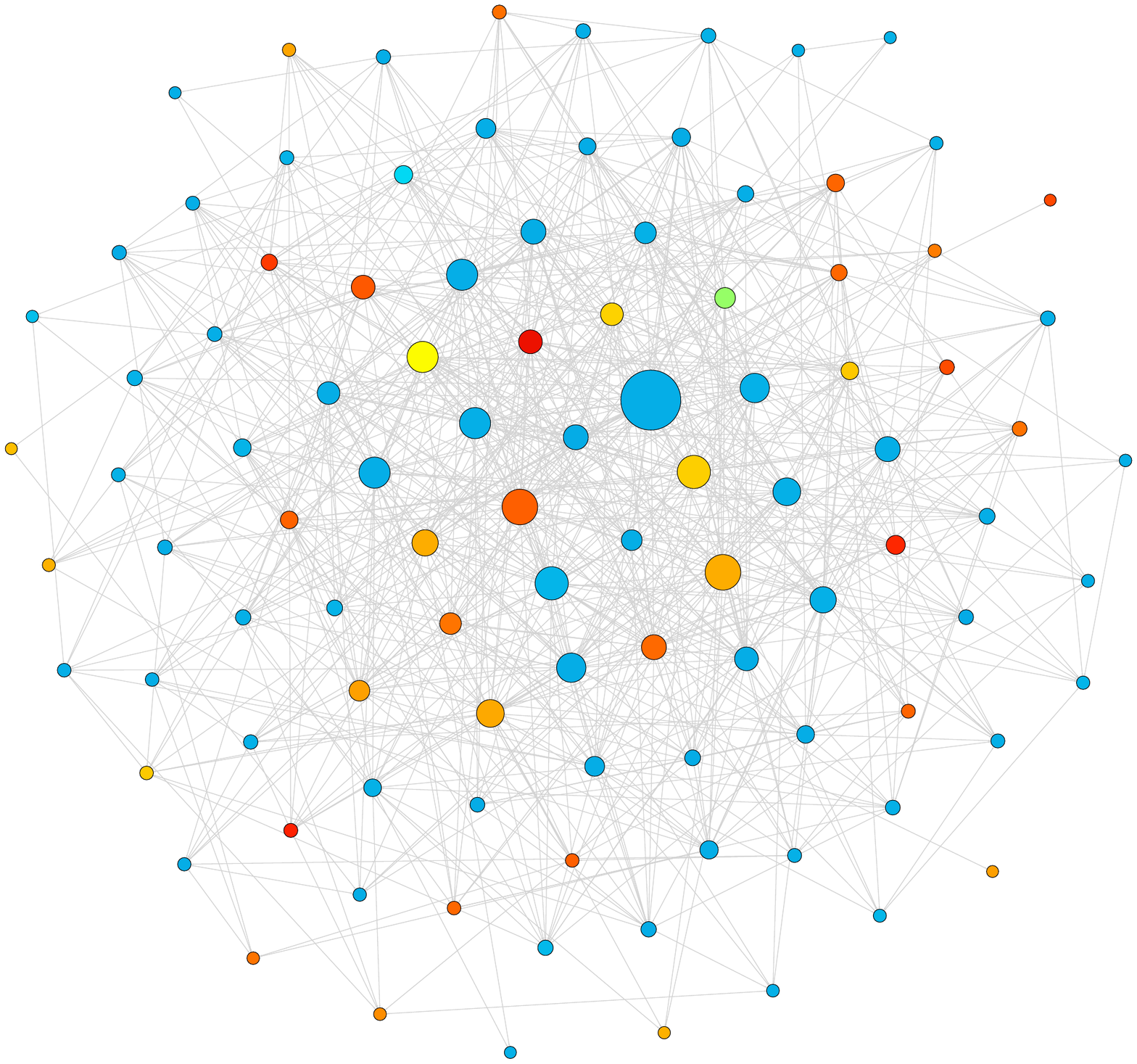}
\includegraphics[scale=0.18]{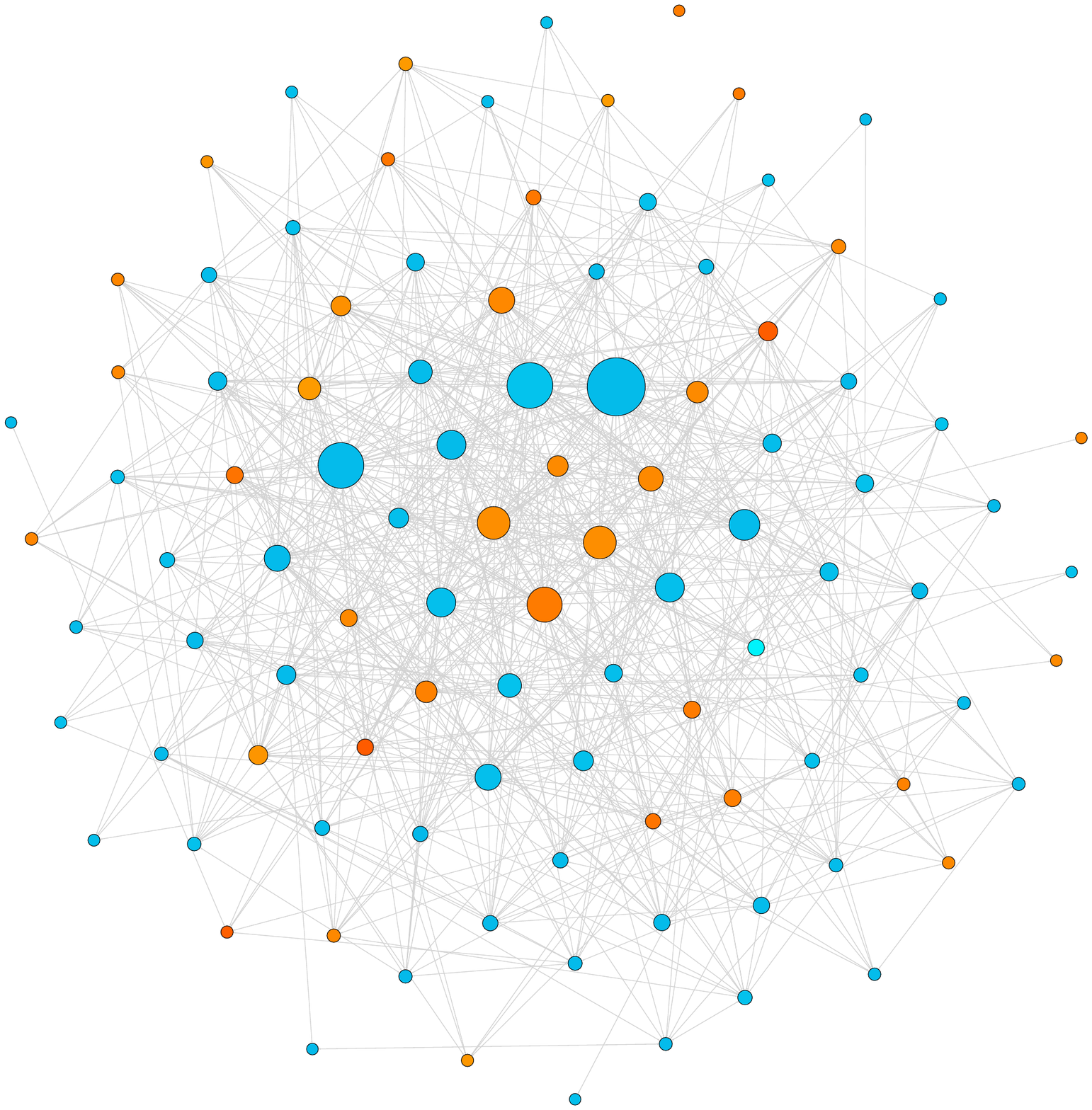} \\
\caption{Evolution of opinion and connection degree of each node of the previously evolved graph $\mathcal{G}^{100}$. From left to right: graph at times $t=0,25,50$. From the top: opinion dynamics for threshold values $c^*=10,20,30$. The target opinion is set $w_d=0.8$ and the control parameter $\kappa=0.1$. } \label{fig:opinion_evo_3}
\end{figure}

\begin{figure}
\centering
\includegraphics[scale=0.36]{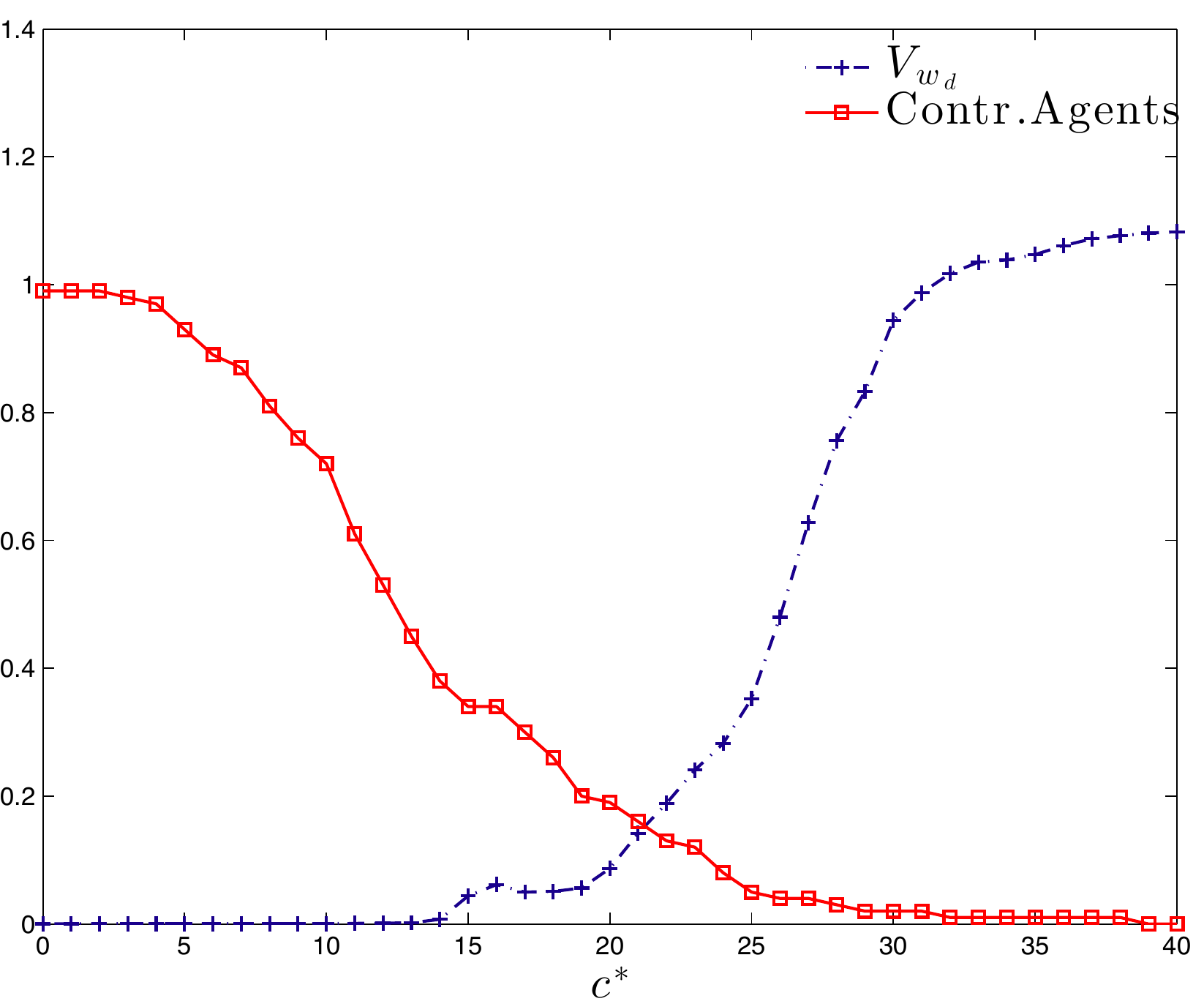}
\includegraphics[scale=0.36]{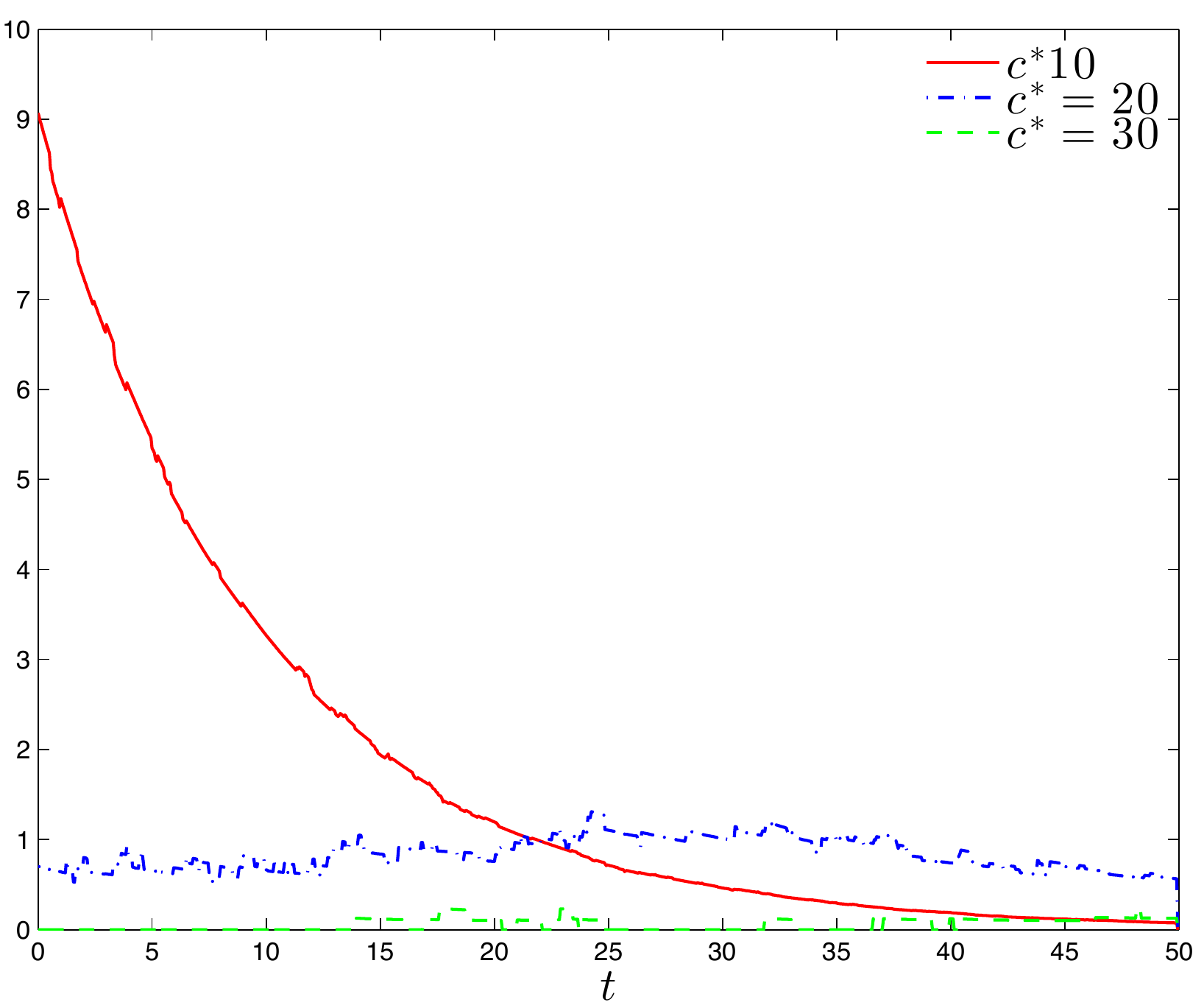}
\caption{Left: the red squared plot indicates the size of the set of controlled agent at the final time $t_f$ in dependence on $c^*$ whereas the blue line indicates the mean square displacement $V_{w_d}$. Right: values of the control $u$ at each time step for $c^*=10,20,30$.  In the numerical test we assumed $\Delta=0.4, \Delta t=5~10^{-3}, \kappa=0.1$.}\label{fig:bounds_cstar}
\end{figure}

In this section we present some numerical results in order to show the main features of the control introduced in the previous paragraphs. We considered a population of $N=100$ agents, each of them representing a node of an undirected graph with density of connectivity $\gamma=30$. The network $\mathcal{G}^{100}$ evolves in the time interval $[0,50]$ with attraction coefficient $\alpha=0.01$ and represents a single sample of the evolution of the master equation \eqref{eq:master} with $D=20$. The control problem is solved by the instantaneous control method described in Remark \ref{rem:inst} with $\Delta t=5~10^{-2}$.  
In Figure \ref{fig:opinion_evo_3} we present the evolution over the reference time interval of the constrained opinion dynamics. The interaction terms have been chosen as follows
\begin{equation}\label{eq:KH}
K(c_i,c_j)=e^{-\lambda c_i}\left(1-e^{-\beta c_j}\right), \qquad H(w_i,w_j)=\chi(|w_i-w_j|\le \Delta),
\end{equation}
where the function $H(\cdot,\cdot)$ is a bounded confidence function with $\Delta=0.4$, while $K(\cdot,\cdot)$ defines the interactions between the agents $i$ and $j$ taking into account that agents with a large number of connections are more difficult to influence and at the same time they have more influence over other agents. The action of the control is characterized by a parameter $\kappa=0.1$ and target opinion $w_d=0.8$. We present the resulting opinion dynamics for a choice of constants $\lambda=1/100,\beta=1$ in Figure \ref{fig:line_nonlin2}. We report the evolution of the network and of the opinion in Figure \ref{fig:opinion_evo_3}, here the diameter of each node is proportional with its degree of connection whereas the color indicates its opinion. 
As a measure of consensus over the agents we introduce the quantity
\begin{equation}
V_{w_d}=\dfrac{1}{N-1}\sum_{i=1}^N (w_i(t_f)-w_d)^2,
\end{equation}
where $w_i(t_f)$ is the opinion of the $i$th agent at the final time $t_f$. In Figure \ref{fig:bounds_cstar} we compare different values of $V_{w_d}$ as a function of $c^*$. Here we calculated the size of the controlled agents and the values of $V_{w_d}$ both, starting from a given uniform initial opinion and the same graph with initial uniform degree distribution.
It can be observed how the control is capable to drive the overall dynamics  toward the desired state acting only on a portion of the nodes.

\section*{Conclusions and perspectives}
In this short note we focus our attention on a control problem for the dynamic of opinion over a time evolving network. We show that the introduction of a suitable selective control depending on the connection degree of the agent's node is capable to drive the overall opinion toward consensus. In a related work we will consider this problem in a mean-field setting where the number of agents, and therefore nodes, is very large \cite{APZc}. 

\section*{Acknowledgments}  
GA acknowledges the support of the ERC-Starting Grant project High-Dimensional Sparse Optimal Control (HDSPCONTR).

\end{document}